\documentclass[12pt]{article}

\usepackage{amssymb}
\usepackage{amsmath}
\usepackage{graphics}
\usepackage{epsfig}
\def\d1{\displaystyle}

\def\t1{\partial}

\def\th(\theta)
\def\e{\varepsilon}

\begin{document}
\begin{center}
\renewcommand{\thefootnote}{\fnsymbol{footnote}}
 {\large\bf  On Liouville Type of Theorems  to the 3-D
 }\\ [3mm]
{\large\bf Incompressible Axisymmetric Navier-Stokes Equations}

\end{center}
\begin{center}

{\normalsize{Quansen JIU }\footnote{The research is partially
supported by National Natural Sciences Foundation of China (No.
11171229, No.11231006 and No.11228102) and Project of Beijing Chang
Cheng Xue Zhe. e-mail:
 jiuqs@mail.cnu.edu.cn}}

{\normalsize\it School of Mathematical Sciences,

Capital Normal University, Beijing 100048,P. R. China}

{\normalsize{ Zhouping XIN }\footnote{The research is partially
supported by Zheng Ge Ru Funds, Hong Kong RGC Earmarked Research
Grants CUHK4041/11P and CUHK4048/13P,  NSFC/RGC Joint Research Scheme Grant
CUHK443/14, a Focus Area Grant at The Chinese University of Hong Kong, and a grant from the Croucher Foundation. e-mail: zpxin@ims.cuhk.edu.hk}}

{\normalsize\it IMS and Department of Mathematics,

The Chinese University of Hong Kong, Shatin, N.T., Hong Kong }

\end{center}

\begin{center}
\parbox{0.8\hsize}{\parindent=4mm
  {\it Abstract}:\  Liouville type of theorems play  a key role in the blow-up approach to study the global regularity of the three-dimensional Navier-Stokes equations.  In this paper, we will prove   Liouville type  of theorems to the 3-D axisymmetric Navier-Stokes equations
with swirls under some   suitable assumptions on swirl component velocity $u_\theta$ which are scaling invariant. It is known that $ru_\theta$ satisfies the maximum principle. The assumptions on $u_\theta$ will be natural and useful to make further studies on the global regularity to the three-dimensional incompressible axisymmetric Navier-Stokes equations.


 {\it AMS(1991)Subject Classification}: 35Q35.}
\end{center}

\section{Introduction}
\setcounter{equation}{0} We consider the  Cauchy problem for the
three-dimensional (3-D) incompressible Navier-Stokes equations
\begin{equation}
\left\{
\begin{array}{ll}\label{1.1}
&\partial_tu-\Delta u+(u\cdot\nabla)u+\nabla p=0, \quad (x,t)\in
\mathbb{R}^3\times (0,T),\\ [3mm] &{\rm div}\  u=0 ,
\end{array}
\right.
\end{equation}
with the initial conditions
\begin{equation}\label{1.3}
u(x,t)\mid_{t=0}=u_0(x).
\end{equation}
Here the unknown functions  are the velocity vector $u=(u_1(x,t), u_2(x,t),
  u_3(x,t))$ and the pressure  $p(x,t)$ with $x\in \mathbb{R}^3, t\in [0,T]$, where $T>0$ is a constant. In (\ref{1.1}),
   $ {\rm div}\  u=0$ means that the fluid is incompressible.

The global existence of the Leray-Hopf weak
solutions to the Cauchy problem or the initial-boundary problem of the three-dimensional Navier-Stokes equations  has been proved long ago (see \cite{Le},\cite{Ho}). However, the uniqueness and
regularity of the weak solutions remain completely open. Up to now,    the weak solutions
will be regular and unique provided that the Serrin-type conditions
  $u\in L^p([0,T);L^q(\mathbb{R}^3))$ hold, where
 $2/p+3/q\le 1$, $p\ge 2$ and $ q\le 3$ (see \cite{Se}, \cite{Str}, \cite{ESS}). The  strong (or smooth) solution of the three-dimensional Navier-Stokes equations was proved to be unique but
 local in time (see \cite{Kat,KT,La,MB,Te}). On the other hand, Scheffer \cite{Sch} introduced and began to study the partial regularity of    suitable weak solutions.  The significant results, due to  Caffarelli-Kohn-Nirenberg
   \cite{CKN} show that, for any  suitable weak solutions, one-dimensional Hausdorff measure  of the singular set is zero. The simplified proofs and further studies   are referred to \cite{Lin}, \cite{TX}.

For the three-dimensional axisymmetric Navier-Stokes equations, if the angular component of the velocity $u_\theta= 0$,  the global existence and uniqueness of the strong (or smooth) solution have been successfully obtained (\cite{La}, \cite{UY}). In the presence of  swirls, that is, $u_\theta\not\equiv 0$, the global well-posedness of the  solution  is still open. Recently, using DeGeogi-Nash-Moser iterations and a blow-up approach respectively,  Chen-Strain-Tsai-Yau \cite{CSTY1,CSTY2} and Koch-Nadirashvili-Seregin-$\breve{S}$ver$\acute{a}k$ \cite{KNSS} obtained an interesting and important development on this problem. Roughly speaking, they proved that if the solution satisfies $(1) \ |ru(x,t)|\le C$ or $(2)\ |u(x,t)|\le \frac{C}{\sqrt{T^*-t}}$ for $0<t<T^*$, where $C>0$ is an arbitrary and absolute constant and $(0,T^*)$ is the maximal existence interval of the solution, then there exists a constant $M>0$ such that $|u(x,t)|\le M$ for $0<t\le T^*$ which implies that the solution is globally regular on time. It should be remarked that these conditions are scaling invariant and imply the possible blow-up rate of the solution. Moreover, the singularity  satisfying (2) is usually called type I singularity in the sense of \cite{Ham}. Thus,   if an axisymmetric solution develops a singularity, it can only be a singularity of type II (any singularity which is  not type I). The other regularity criteria and recent studies  can be seen in \cite{Chae, JX, KPZ, LZ, Pan} and references therein.

 The basic idea of the blow-up approach is that  if the solution would blow up at some space-time point, then making scaling transformation of the solution and enlarging the region near the possibly singular  point, one can  look into the equations satisfied by some suitably  scaled  solutions. In particular, after taking the limit, if the  solution of the limit equation is trivial, which is so called a Liouville type of theorem, then one will obtain a contradiction and the blow-up will not happen. To the three-dimensional
 axisymmetric Navier-Stokes equations, the possible singularity of the solution may only appear on the symmetry axis due to the partial regularity theory in \cite{CKN}. Therefore it suffices to study the possible singularity of the solution on the symmetry axis. In this paper, we are concerned with   Liouville type of theorems.

First,  we prove a Liouville type of theorem    by assuming that
\begin{eqnarray}\label{Oct1}
\displaystyle\limsup_{\delta\to
0+}\|ru_\theta(x,t)\|_{L^\infty(\{x|r\le
\delta\}\times (-1,0))}=0.
\end{eqnarray}
It is shown that, under the assumption  \eqref{Oct1}, there exists a bounded and continuous function $s(T)$ defined on $(-\infty,0]$ such that the ancient solution $\bar u=(0,0, s(T))$.  It should be remarked that the assumption of \eqref{Oct1} is  natural since $ru_\theta$ satisfies the maximum principle and  if the initial data satisfies $|ru_{\theta0}|\le C$ then the solution will keep the bound $|ru_\theta|\le C$ for some constant $C>0$. This implies that the  singularity of $u_\theta$ near the symmetry axis, if exists,  may be  of the  rate $\frac{O(1)}{r}$ as $r\to 0$, where $O(1)$ means a finite constant. While  the condition \eqref{Oct1} implies that the  singularity which we impose on $u_\theta$ near the symmetry axis is  of the rate $\frac{o(1)}{r}$ with $o(1)\to 0$ as $r\to 0$.    Our approaches are based on  \cite{CSTY1,CSTY2} and  \cite{KNSS}. In particular, we will use the integral expression on $\frac{\omega_\theta^{(k)}}{R}$ which is the scaled quantity of $\frac{\omega_\theta}{r}$ to prove a Liouville type of theorem  to the ancient solution, where $\omega_\theta=\partial_ru_3-\partial_3u_r$ is the angular component of the vorticity.  This is different from \cite{KNSS} in which the authors established a Liouville type of theorem by making full use of the strong maximum principle of the scalar equation of  $\Gamma=ru_\theta$, under the assumption that $|ru|\le C$. Moereover, in comparison with the global regularity results in \cite{CSTY1,CSTY2} and  \cite{KNSS}, we need \eqref{Oct1} but do not require the condition on the radial component of the velocity $u_r$. It should be noted that in the process of proving the Liouville type of theorem, we only need the condition \eqref{Oct1}.  How to remove the condition \eqref{Oct1}  will be very interesting and challenging.

Second,  we  prove a Liouville type of theorem  under weighted estimates of smooth solutions to the three-dimensional
 axisymmetric Navier-Stokes equations, which are also scaling invariant. This is motivated by weighted estimates in \cite{CKN} which can be carried out under assumptions of the initial data. Finally, in the end of the paper, we give further remarks on the possibility to rule out the singularity of the solution by using the Liouvile type of theorems established in this paper.

The organization of this paper is as follows. In Section 2, we will present some
preliminaries and establish a Liouville type of theorem under \eqref{Oct1}. In Section 3,  we  prove a Liouville type of theorem  under weighted estimates of the  solution. In Section 4,  we give further remarks  on how to use the Liouville type of theorem in future works.

\section{A Liouville Type of Theorem Under \eqref{Oct1}}
\setcounter{equation}{0}

By an axisymmetric solution $(u, p)$ of \eqref{1.1}, we mean that,
 in the cylindrical coordinate systems, the solution takes the form  $p(x,t)=p(r,x_3,t)$ and
$$
u(x,t)=u_r(r,x_3,t)e_r+u_\theta(r,x_3,t)e_\theta+u_3(r,x_3,t)e_3,
$$
where
$$
e_r=(\cos\theta, \sin\theta, 0),\quad  e_\theta=(-\sin\theta,
\cos\theta, 0),\quad e_3=(0, 0, 1).
$$
Here $u_\theta(r,x_3,t)$ and  $u_r(r,x_3,t)$ are  the angular and radial components of $u(x,t)$ respectively. For
the axisymmetric velocity field $u$, the corresponding vorticity
$\omega=\nabla\times u$ is
$$
\omega=\omega_r e_r+\omega_\theta e_\theta+\omega_3 e_3,
$$
where
$$
\omega_r=\t1_3 u_\theta,\  \omega_\theta=\t1_r u_3-\t1_3 u_r, \
\omega_3=-\frac 1r\t1_r(ru_\theta).
$$
The 3-D axisymmetric  Navier-Stokes equations read as
\begin{eqnarray}
 &&\frac{\tilde Du_r
}{Dt}-(\partial_r^2+\partial_3^2+\frac{1}{r}\partial_r)u_r+\frac{1}{r^2}u_r-\frac{1}{r}(u_\theta)^2+\partial_rp=0,
\label{1.4} \\[3mm]
&&\frac{\tilde
Du_\theta}{Dt}-(\partial_r^2+\partial_3^2+\frac{1}{r}\partial_r)u_\theta+\frac{1}{r^2}u_\theta
+\frac{1}{r}u_\theta u_r=0, \label{1.5} \\[3mm]
&&\frac{\tilde
Du_3}{Dt}-(\partial_r^2+\partial_3^2+\frac{1}{r}\partial_r)u_3+\t1_3p=0,
\label{1.6} \\[3mm]
&&\t1_r(ru_r)+\t1_3(ru_3)=0, \label{1.7}
\end{eqnarray}
where
$$
\frac{\tilde D}{Dt}=\t1_t+u_r\t1_r+u_3\t1_3, \quad
r=(x_1^2+x_2^2)^{1/2}.
$$

In the following, we set
$$
\tilde\nabla=(\partial_r, \partial_3)
$$
and use $C$ to denote an absolute constant
which may be different from line to line.

Without loss of generality, after translation on the time variable, $u(x,t)$ is assumed to be a smooth axisymmetric solution to \eqref{1.4}-\eqref{1.7},
defined in   $\mathbb{R}^3\times (-1,0)$ with $u\in L^\infty(\mathbb{R}^3\times (-1,t')$ for any $-1<t'<0$.

Let
\begin{eqnarray}\label{P1}
h(t)=\sup_{x\in \mathbb{R}^3}|u(x,t)|, \ \ H(t)=\sup_{-1\le s\le t<0} h(s).
\end{eqnarray}

 Suppose that the first singularity time for the solution $u(x,t)$ is at time $t=0$.
 Then it is clear that $\lim_{t\to 0-}
H(t)=\infty$. In fact, by a classical result of Leray \cite{Le},
if $u$ develops a singularity at $t=0$, then
\begin{eqnarray}\label{P2}
h(t)=\sup_{x\in \mathbb{R}^3}|u(x,t)|\ge \frac{\varepsilon_1}{\sqrt{-t}}
\end{eqnarray}
for some $\varepsilon_1>0$.

There exist $t_k\nearrow 0$ as $k\to\infty$ such that
$H(t_k)=h(t_k)$. Denote $N_k=H(t_k)$. Then there exists a sequence
of numbers $\gamma_k\searrow 1$ as $k\to\infty$ and $x_k\in \mathbb{R}^3$
such that $M_k=|u(x_k,t_k)|\ge N_k/\gamma_k, k=1,2\cdots$,
satisfying $M_k\to \infty$ as $k\to\infty$.

Define
\begin{eqnarray}\label{P5}
u^{(k)}(X,T)=\frac{1}{M_k}u(\frac{X_1}{M_k}, \frac{X_2}{M_k},
x_{k3}+\frac{X_3}{M_k}, t_k+\frac{T}{M^2_k}), k=1,2,\cdots
\end{eqnarray}
In the cylindrical coordinate system, set
$$
u^{(k)}(X,T)=b^{(k)}(X,T)+u_\theta^{(k)}e_\theta,
$$
where $b^{(k)}(X,T)=u^{(k)}_Re_R+u^{(k)}_Ze_Z, R=\sqrt{X_1^2+X_2^2}$.

 Then $u^{(k)}(X,T)$ are  smooth solutions of the 3D Navier-Stokes equations, which are defined in
$\mathbb{R}^3\times (A_k,B_k)$ with
\begin{eqnarray}\label{P60}
A_k=-M_k^2-M_k^2t_k, \ B_k=-M_k^2t_k.
\end{eqnarray}
Note that $B_k=-M_k^2t_k\ge (\frac{N_k}{\gamma_k})^2(-t_k)\ge
\frac{\varepsilon_1}{\gamma_k^2}$. Moreover, it holds that
\begin{eqnarray}\label{P7}
|u^{(k)}(X,T)|\le \gamma_k, X\in \mathbb{R}^3, T\in (A_k,0),
\end{eqnarray}
and
\begin{eqnarray}\label{P8}
 |u^{(k)}(M_kx_{k1},M_kx_{k2},0,0)|=1.
\end{eqnarray}
It follows from the regularity theorem of the Navier-Stokes
equations that
\begin{eqnarray}\label{P8+}
 |\partial_Tu^{(k)}|+|D^lu^{(k)}|\le C_l, X\in \mathbb{R}^3, T\in (A_k,0]
\end{eqnarray}
for $k=1, 2,\cdots$ and $|l|=0, 1, 2,\cdots$, where $l=(l_1,l_2,l_3)$ is a multi-index satisfying $l_1+l_2+l_3=|l|$ and $D^l=\frac{\partial^{|l|}}{\partial x_1^{l_1}\partial x_2^{l_2}\partial x_3^{l_3}}$.  $C_l$ is a constant depending on $l$ but not on $k$.
Then there exists a smooth function $\bar u(X,T)$ defined in $\mathbb{R}^3\times (-\infty,0)$ such that, for any $|l|=0,1,2 \cdots$,
\begin{eqnarray}\label{Oct26}
D^l u^{(k)}\longrightarrow D^l \bar u, \ k\to \infty, \
\end{eqnarray}
uniformly in $C(\bar Q)$ for any compact subset $Q\subset\subset
\mathbb{R}^3\times (-\infty,0]$.  Denote $\bar \omega(X,T)=\bar\omega_\theta e_\theta+\bar\omega_re_r+\bar\omega_ze_z$  the voricity of $\bar u(X,T)$.

\vspace{2mm}

Our   main result of this section is a Liouville type of  theorem as follows.

{\bf Theorem 2.1} Let $u(x,t)$ be an axisymmetric vector field
defined in $\mathbb{R}^3\times(-1,0)$ which belongs to
$L^\infty(\mathbb{R}^3\times(-1,t')$ for each $-1<t'<0$. Assume that $u$ satisfies
 \begin{eqnarray}\label{C10-1}
|ru_\theta(x,t)|\le C, \ \ (x,t)\in \mathbb{R}^3\times (-1,0),
\end{eqnarray}
and
\begin{eqnarray}\label{Oct30}
&\displaystyle\limsup_{\delta\to
0+}\|ru_\theta(x,t)\|_{L^\infty(\{x|r\le
\delta\}\times (-1,0))}=0,\label{C1+}
\end{eqnarray}
where $C>0$ is any finite constant. Then either
\begin{eqnarray}\label{A-10}
|u(x,t)|\le M, \ \ x\in \mathbb{R}^3,\ \  t\in [-1,0],
\end{eqnarray}
where $M>0$ is an absolute constant depending on $C$, or $\bar\omega=0$ and
$\bar u=(0,0,s(T))$, where $\bar u$ is same as in \eqref{Oct26} and $s(T): (-\infty, 0]\to \mathbb{R}$ is a bounded and continuous function.

{\bf Remark 2.1} The condition (\ref{C10-1})  can be removed if the initial data satisfies
$\|ru_{0\theta}\|_{L^\infty}<\infty$ (see \cite{Chae,JX}). The condition \eqref{Oct30} means that the  singularity of $u_\theta$ near the symmetry axis is of the rate $\frac{o(1)}{r}$ with $o(1)\to 0$ as $r\to 0$.

{\bf Proof.} Suppose that \eqref{A-10} is false. Then one can rescale the solution as in \eqref{P5}-\eqref{Oct26}. It will be shown that $\bar\omega(X,T)=0$ and
$\bar u(X,T)=(0,0,s(T))$ with $s(T): (-\infty, 0)\to \mathbb{R}$  a bounded and continuous function.

Let $C_0>0$ be any fixed constant. For any $X\in \mathbb{R}^3$ with $R\le C_0$ and $T\in (A_k,0]$, it follows from
\eqref{C1+} that
\begin{eqnarray}\label{P14}
&&|\Gamma^{(k)}(X,T)|\equiv|Ru_\theta^{(k)}|=|\frac{R}{M_k}u_\theta(\frac{X_1}{M_k},
\frac{X_2}{M_k}, x_{k3}+\frac{X_3}{M_k},
t_k+\frac{T}{M^2_k})|\nonumber\\[3mm]
&&\le \tilde F(k,C_0)\to 0
\end{eqnarray}
as $k\to \infty$.  Set
$$
F(k,C_0)=\max(\tilde F(k,C_0),\frac1k), \quad k=1,2,3\cdots.
$$
It follows that
\begin{eqnarray}\label{P15}
|u_\theta^{(k)}(X,T)|\le R^{-1}F(k,C_0), 0<R\le C_0, T\in (A_k,0].
\end{eqnarray}

It follows from (\ref{C10-1}) that
$$
|u_\theta^{(k)}(X,T)|\le \frac{C}{R}, \ \ R>0,
 T\in (A_k,0]
$$
Using (\ref{P8+}) and  the fact that $u_\theta^{(k)}|_{R=0}=0$, one has
\begin{eqnarray}
&|u_\theta^{(k)}(X,T)|\le C \min(R, R^{-1}), R>0, T\in
(A_k,0],\label{P16}\\[3mm]
& |\partial_Zu_\theta^{(k)}(X,T)|\le C \min(R, 1), R>0, T\in
(A_k,0].\label{P17}
\end{eqnarray}
Consequently, for any $T\in (A_k,0]$,
\begin{eqnarray}\label{June-2-1}
|u_\theta^{(k)}(X,T)|\le \left\{
\begin{array}{lll}
 &CR, \ & R<\sqrt{F(k,C_0)}, \\[3mm]
&\frac{F(k,C_0)}{R}, \ & \sqrt{F(k,C_0)}\le R\le C_0, \\[3mm]
&\frac{C}{R}, \ &R>C_0.
\end{array}
 \right.
\end{eqnarray}
\begin{eqnarray}\label{June-2-2}
|\frac{\partial_Z(u_\theta^{(k)})^2(R,Z,T)}{R^2}|\le \left\{
\begin{array}{lll}
 &C, \ & R<\sqrt{F(k,C_0)}, \\[3mm]
&C\frac{F(k,C_0)}{R^2}, \ &\sqrt{F(k,C_0)}\le R<1,\\[3mm]
&C\frac{F(k,C_0)}{{R}^3}, \ & 1\le R\le C_0, \\[3mm]
&\frac{C}{{R}^3}, \ &R\ge C_0.
\end{array}
\right.
\end{eqnarray}

Let $\Omega=\frac{\omega_\theta(x,t)}{r}$ and
$f^{(k)}=\Omega^{(k)}(X,T)=\frac{\omega^{(k)}_\theta(X,T)}{R}$. Then
it holds that
\begin{eqnarray}\label{20}
|f^{(k)}(X,T)|\le C(1+R)^{-1},\ \  X\in \mathbb{R}^3, T\in
(A_k,0).
\end{eqnarray}
It follows from the equation of $\omega_\theta$ that
\begin{eqnarray}\label{21}
(\partial_T-L)f^{(k)}=g^{(k)}, \ \
L=\Delta+\frac2R\partial_R-b^{(k)}\cdot\nabla_X,
\end{eqnarray}
where $g^{(k)}=R^{-2}\partial_Z(u_\theta^{(k)})^2$ and
$b^{(k)}=u_R^{(k)}e_R+u_Z^{(k)}e_Z$.

Regarding $f^{(k)}(X,T)=f^{(k)}(R,Z,T)$ as a 5-dimensional
axisymmetric function by denoting $X=(\tilde X, X_5)=(X_1,\cdots,X_4,X_5)$,  $R=|\tilde X|=\sqrt{X_1^2+X_2^2+X_3^2+X_4^2}$ and $Z=X_5$,
we obtain
\begin{eqnarray}\label{21+}
(\partial_T+\tilde b^{(k)}\cdot\tilde
\nabla_X-\Delta_5)f^{(k)}=g^{(k)}, \ \
\end{eqnarray}
where $\tilde b^{(k)}=u_R^{(k)} \tilde e_R+u_Z^{(k)}\tilde e_Z$ with
$\tilde e_R=(\frac{X_1}{R}, \frac{X_2}{R}, \frac{X_3}{R},
\frac{X_4}{R}, 0)$ and $\tilde e_Z=(0,0,0,0,1)$. The scaling \eqref{P5} can be rewritten as
\begin{eqnarray}\label{Jan2}
u_\theta^{(k)}(Y,T)=\frac{1}{M_k}u(\frac{R}{M_k},
z_k+\frac{Y_5}{M_k},t_k+\frac{T}{M_k^2}),
\end{eqnarray}
where $z_k=x_{3k}$.

Denote $P(T,X;S,Y)$ the kernel for $\partial_T+\tilde
b^{(k)}\cdot\tilde \nabla_X-\Delta_5$ and $Y=(\tilde Y, Y_5)=(Y_1, \dots, Y_5)$. By the Duhamel's formula,
\begin{eqnarray}\label{22}
&&f^{(k)}(X,T)=\int P(T,X;S,Y)f^{(k)}(Y,S) dY+\int_S^T\int
P(T,X;\tau,Y)g^{(k)}(Y,\tau) dYd\tau\nonumber\\[3mm]
&&=:I+II.
\end{eqnarray}

Due to Carlen-Loss \cite{CL} and Chen-Strain-Tsai-Yau \cite{CSTY2}, the kernel $P$ satisfies $P\ge 0, \int P(T,X;S,Y) dY=1$ and
$$
P(T,X;S,Y)\le C(T-S)^{-\frac52}e^{-h(|X-Y|,T-S)}, \
h(a,T)=C\frac{a^2}{T}[(1-\frac{T}{a})_+]^2,
$$
where $f_+=\max \{0, f\}$ and we have used the fact that $\|\tilde b^{(k)}\|_\infty\le \gamma_k\le 2$ for $k\ge N$.

It can be verified that the function $e^{-h(a,T)}, a\ge 0, T\ge 0,$ has the following properties:

When $T\ge T_0>0$ for any fixed (but may be small) $T_0>0$, one has
\begin{eqnarray}\label{23-0}
e^{-h(a,T)}\le Ce^{-Ca/T}
\end{eqnarray}
holds for some constant $C>0$ which may depend on $T_0$.

When $T\ge \frac a2$, one has
\begin{eqnarray}\label{23-1}
e^{-h(a,T)}\le Ce^{-Ca/T}.
\end{eqnarray}

When $0\le T\le \frac a2$, it is easy to get
\begin{eqnarray}\label{23-2}
e^{-h(a,T)}=e^{-C\frac{a^2}{T}(1-\frac{T}{a})^2}\le e^{-Ca^2/T}.
\end{eqnarray}
In \eqref{23-1} and \eqref{23-2}, $C>0$ is some uniform constant.

It follows from (\ref{23-0}) and H\"older inequality that
\begin{eqnarray}\label{24}
&&|I|\le [\int P(T,X;S,Y)|f^{(k)}(Y,S)|^5
dY]^{\frac15}\nonumber\\[3mm]
&&\le [C(T-S)^{-\frac52}\int e^{-C\frac{|X_5-Y_5|}{T-S}}\frac{\mathbb{R}^3
dR}{(1+R)^5}
dY_5]^\frac15\nonumber\\[3mm]
&&\le C(T-S)^{-\frac{3}{10}}
\end{eqnarray}
for all $X\in \mathbb{R}^5, T, S\in (A_k,0)$ satisfying $T-S>0$.

Let $L=\{\tau\in [S,T]: T-\tau\ge \frac{|X-Y|}{2}\}$ and
$L^c=\{\tau\in [S,T]: T-\tau\le \frac{|X-Y|}{2}\}$.  For any
$C_0>0$, with help of (\ref{23-1}) and (\ref{23-2}), we have
\begin{eqnarray*}\label{25}
&&|II|\le \int_L\int C(T-\tau)^{-\frac52}
e^{-C\frac{|X-Y|}{T-\tau}}|g^{(k)}(Y,\tau)|
dYd\tau\nonumber\\[3mm]
&&+\int_{L^c} \int C(T-\tau)^{-\frac52}
e^{-C\frac{|X-Y|^2}{T-\tau}}|g^{(k)}(Y,\tau)|
dYd\tau\nonumber\\[3mm]
 &&\le C\int_S^T (T-\tau)^{-\frac52}
d\tau\int e^{-C\frac{|X_5-Y_5|}{T-\tau}}dY_5(\int_{\{R\le
C_0\}}e^{-C\frac{|\tilde X-\tilde Y|}{T-\tau}} d\tilde Y)^{\frac
12}(\int_{\{ R\le C_0\}} |g^{(k)}(Y,\tau)|^2 d\tilde Y)^{\frac
12}\nonumber\\[3mm]
&& +C\int_S^T (T-\tau)^{-\frac52} d\tau\int
e^{-C\frac{|X_5-Y_5|}{T-\tau}}dY_5(\int_{\{R\ge
C_0\}}e^{-C\frac{|\tilde X-\tilde Y|}{T-\tau}} d\tilde Y)^{\frac
12}(\int_{\{ R\ge C_0\}} |g^{(k)}(Y,\tau)|^2 d\tilde Y)^{\frac 12}\nonumber\\[3mm]
&&+C\int_S^T (T-\tau)^{-\frac52}d\tau\int
e^{-C\frac{|X_5-Y_5|^2}{T-\tau}} dY_5(\int_{\{R\le C_0\}}
e^{-C\frac{|\tilde X-\tilde Y|^2}{T-\tau}} d\tilde
Y)^{\frac23}(\int_{\{ R\le C_0\}} |g^{(k)}(Y,\tau)|^3 d\tilde
Y)^{\frac
13}\nonumber\\[3mm]
&&+C\int_S^T (T-\tau)^{-\frac52}d\tau\int
e^{-C\frac{|X_5-Y_5|^2}{T-\tau}} dY_5(\int_{\{R\ge C_0\}}
e^{-C\frac{|\tilde X-\tilde Y|^2}{T-\tau}} d\tilde
Y)^{\frac23}(\int_{\{ R\ge C_0\}} |g^{(k)}(Y,\tau)|^3 d\tilde
Y)^{\frac
13}.\nonumber\\[3mm]
\end{eqnarray*}
It follows from (\ref{June-2-2}) that
\begin{eqnarray*}
&&\int_{\{ R\le C_0\}} |g^{(k)}(Y,\tau)|^2 d\tilde Y\\[3mm]
&&=[\int_{\{ R<\sqrt{F(k,C_0)}\}}+\int_{\{ \sqrt{F(k,C_0)}\le
R<1\}}+\int_{\{
1\le R<C_0\}}] |g^{(k)}(Y,\tau)|^2 d\tilde Y\\[3mm]
&&\le C[(\sqrt{F(k,C_0)})^4+F^2(k,C_0)(-\ln
\sqrt{F(k,C_0)})+F^2(k,C_0)].
\end{eqnarray*}
and
\begin{eqnarray*}
&&\int_{\{ R\le C_0\}} |g^{(k)}(Y,\tau)|^3 d\tilde Y\\[3mm]
&& \le C[(\sqrt{F(k,C_0)})^4+F^3(k,C_0)(
\sqrt{F(k,C_0)})^{-2}+F^3(k,C_0)].
\end{eqnarray*}
It concludes that
\begin{eqnarray*}
|II|&&\le
C(T-S)^\frac32[(\sqrt{F(k,C_0)})^4+F^2(k,C_0)(-\ln
\sqrt{F(k,C_0)})+F^2(k,C_0)+C_0^{-1}]\nonumber\\[3mm]
&&+C(T-S)^\frac 13[(\sqrt{F(k,C_0)})^4+F^3(k,C_0)(
\sqrt{F(k,C_0)})^{-2}+F^3(k,C_0)+C_0^{-\frac 53}].
\end{eqnarray*}
Letting $T-S=C_0^\frac13$, we obtain
\begin{eqnarray*}\label{25+}
&& |I|+|II|\le
C(C_0^{-\frac{1}{10}}+C_0^{-\frac12}+C_0^{-\frac{14}{9}})\\[3mm]
&&+C[(\sqrt{F(k,C_0)})^4+F^2(k,C_0)(-\ln
\sqrt{F(k,C_0)})+F^2(k,C_0)]^\frac12\nonumber\\[3mm]
&&+C[(\sqrt{F(k,C_0)})^4+F^3(k,C_0)(
\sqrt{F(k,C_0)})^{-2}+F^3(k,C_0)]^\frac13.
\end{eqnarray*}
Taking the limit $k\to\infty$ first and then letting
 $C_0\to \infty$, we have
 \begin{eqnarray*}\label{26+}
|f^{(k)}(X,T)| \to 0,
\end{eqnarray*}
and hence
\begin{eqnarray}\label{26+}
|\omega_\theta^{(k)}(X,T)| \to 0,
\end{eqnarray}
which implies that $\bar\omega_\theta(R,Z,T)=0$. Consequently, we obtain that $\bar \omega_\theta(X,T)=0$ for $X\in \mathbb{R}^3$ and $T\in (-\infty,0)$.
Denote $\bar b(X,T)=\bar u_Re_R+\bar u_Z$. It follows from  \eqref{Oct26} that
\begin{eqnarray}\label{27}
u^{(k)}(X,T)\to  \bar u (X,T),\  b^{(k)}(X,T)\to  \bar
b (X,T), \ {\rm as}\ \ k\to \infty,
\end{eqnarray}
uniformly in $C(\bar Q)$ with any compact subset $Q\subset\subset
\mathbb{R}^3\times (-\infty,0]$.

As a consequence of $curl_X b^{(k)}(X,T)=\omega_\theta^{(k)}e_\theta$,
$div_X b^{(k)}(X,T)=0$ and \eqref{27},  $\bar
b(X,T)$ is a harmonic and bounded function. That is, $\Delta_X \bar b=0$ and $\bar b(X,T)$ is bounded. Since $\bar u_R(0,Z,T)=0$, there exists a continuous and bounded function $s(T): (-\infty,0]\to \mathbb{R}$ such that $\bar b(X,T)=(0, s(T))$. Moreover, since $|Ru_\theta^{(k)}|\to 0$ as $k\to \infty$ for  $0<R\le C_0$, where $C_0>0$ is  arbitrary,  it follows that $\bar u_\theta(X,T)=0$ for all $X\in \mathbb{R}^3, T\in (-\infty, 0)$. Therefore we have proved that $\bar u(X,T)=(0, 0, s(T))$ and furthermore $\bar\omega(X,T)=0$. The proof of the theorem is finished.

\hspace{2mm}

\section{A Liouville Type Theorem Under Weighted Estimates}
\setcounter{equation}{0}

In \cite{CKN}, some weighted estimates were obtained for suitable weak solutions $(u(x,t),p(x,t))$
  of the three-dimensional Navier-Stokes equations \eqref{1.1}.  More precisely, suppose that there exists a small number $L_0>0$
such that if
\begin{eqnarray}\label{Feb22-1}
\int_{\mathbb{R}^3} \frac{|u_0|^2}{|x|} dx \le L
\end{eqnarray}
with $0<L<L_0$, then
\begin{eqnarray}\label{Feb22-2}
\int_{\mathbb{R}^3} \frac{|u|^2}{|x|}
dx+(L_0-L)exp\{\frac{1}{L_0}\int_0^t\int_{\mathbb{R}^3} \frac{|\nabla
u|^2}{|x|}\} dxdt\le L_0, \ \ t\in [0,T].
\end{eqnarray}

For the three-dimensional axisymmetric Navier-Stokes equations, if
\begin{eqnarray}\label{Oct26-3}
\int_{\mathbb{R}^3} \frac{|u_0|^2}{r} dx \le L,
\end{eqnarray}
where $L>0$ is same as in \eqref{Feb22-1}, then it is clear that
 \begin{eqnarray}\label{N-7-1}
\int_{\mathbb{R}^3} \frac{|u_0|^2}{\sqrt{r^2+(z-z_0)^2}} dx \le L
\end{eqnarray}
for any $z_0\in R$. By the translation with respect to $z(=x_3)$,
one can prove in a similar way as in \cite{CKN} that if
\eqref{N-7-1} holds,
 then
\begin{eqnarray}\label{N-7-2}
&&\int_{\mathbb{R}^3} \frac{|u|^2}{\sqrt{r^2+(z-z_0)^2}}
dx\nonumber\\
&&~~~~~~~+(L_0-L)exp\{\frac{1}{L_0}\int_0^t\int_{\mathbb{R}^3}
\frac{|\nabla u|^2}{\sqrt{r^2+(z-z_0)^2}}\} dxdt\le L_0
\end{eqnarray}
for any $z_0\in \mathbb{R}$, where  $L$ and $L_0$ are same as in
\eqref{Feb22-2}.  Note that the quantities on the left hand side of
\eqref{N-7-2} are scaling invariant.

Let $B_{\bar R}=\{x\in \mathbb{R}^3| |x|<R\}$ be a ball with radius $\bar R>0$.
Motivated by the weighted estimates \eqref{N-7-2}, we have

{\bf Theorem 3.1} Let $u(x,t)$ be an axisymmetric vector field
defined in $\mathbb{R}^3\times(-1,0)$ which belongs to
$L^\infty(\mathbb{R}^3\times(-1,t')$ for each $-1<t'<0$ and $L^\infty (B_{\bar R}^c\times (-1,0))$  for some $\bar R>0$, where  $B_{\bar R}^c=\mathbb{R}^3\backslash B_{\bar R}$. Moreover, assume that $u$ satisfies
 \begin{eqnarray}
&& (1)\ \  |ru_\theta(x,t)|\le C, \ \ (x,t)\in \mathbb{R}^3\times (-1,0), \label{C1}\\[3mm]
&& (2)\ \ \int_{\mathbb{R}^3} \frac{|u|^2}{\sqrt{r^2+(z-z_0)^2}}
dx+\int_{-1}^0\int_{\mathbb{R}^3}
\frac{|\omega_r|^2}{\sqrt{r^2+(z-z_0)^2}} r drdzdt\le C \
\label{Oct26-2}
\end{eqnarray}
for any $z_0\in \mathbb{R}$.

 Then either
\begin{eqnarray}\label{A-1}
|u(x,t)|\le M, \ \ x\in \mathbb{R}^3,\ \  t\in (-1,0],
\end{eqnarray}
where $M>0$ is an absolute constant depending on $C$, or
$\bar u=0$, where $\bar u$ is defined in \eqref{Oct26}.

{\bf Remark 3.1} The assumption that $u(x,t)\in L^\infty (B_{\bar R}^c\times (-1,0))$  for some $R>0$ can be easily satisfied if the solution decays at far fields for all $t\in (-1,0)$. In particular, it follows from \cite{CKN} that

{\bf Proposition 3.2}  Suppose that $u_0\in L^2(\mathbb{R}^3)$ and
$$
\int_{\mathbb{R}^3} |\nabla u_0|^2 dx <\infty.$$
Given a  suitable weak solution $u(x,t)$ of the three-dimensional Navier-Stokes equations with initial data $u_0$. Then $u(x,t)$ is regular in the region $B^c_{\bar R}$ for some $\bar R>0$.

\hspace{3mm}

{\bf Remark 3.2} The second term in  \eqref{Oct26-2} can be replaced by
\begin{eqnarray}\label{Oct26-4}
\int_0^t\int_{\mathbb{R}^3} \frac{|\omega_r|^2}{r}
dxdt\le C.
\end{eqnarray}
Note that $\omega_r=\partial_zu_\theta$. The second term in the condition \eqref{Oct26-2} is a weighted estimate on $\partial_zu_\theta$.

\hspace{3mm}
Now we give a proof of Theorem 3.1.

\hspace{3mm}

 {\bf Proof of Theorem 3.1} Suppose that the solution $u(x,t)$ has  singularity at $t=0$. Then  similar to the proof of Theorem 2.1,  there exist $t_k\nearrow 0$ as $k\to\infty$ such that
$H(t_k)=h(t_k)$, where $H(t)$ and $h(t)$ are same as in \eqref{P1}. Denote $N_k=H(t_k)$. Since $u(x,t)\in L^\infty (B_{\bar R}^c\times (-1,0))$  for some $\bar R>0$ by the assumption,  we can choose a sequence
of numbers $\gamma_k\searrow 1$ as $k\to\infty$ and $x_k\in B_{\bar R}$
such that $M_k=|u(x_k,t_k)|\ge N_k/\gamma_k, k=1,2\cdots$,
satisfying $M_k\to \infty$ as $k\to\infty$. Using the same scaling \eqref{P5}, one can prove the theorem as for Theorem 2.1.

To be more precise, we  continue the proof based on  \eqref{22}. The first term $I$ in
\eqref{22} can be estimated as in \eqref{24}. Thus, one needs to focus on the
estimate of the second term $II$ in \eqref{22}.

Note that
\begin{eqnarray}\label{251}
&&|II|\le \int_L\int C(T-\tau)^{-\frac52}
e^{-C\frac{|X-Y|}{T-\tau}}|g^{(k)}(Y,\tau)|
dYd\tau\nonumber\\[3mm]
&&+\int_{L^c} \int C(T-\tau)^{-\frac52}
e^{-C\frac{|X-Y|^2}{T-\tau}}|g^{(k)}(Y,\tau)|
dYd\tau\equiv II_1+II_2,
\end{eqnarray}
where $L=\{\tau\in [S,T]: T-\tau\ge \frac{|X-Y|}{2}\}$ and
$L^c=\{\tau\in [S,T]: T-\tau\le \frac{|X-Y|}{2}\}$. Here we are estimating in 5-dimensional space, $X=(\tilde X,X_5)=(X_1,\cdots,X_4,X_5), Y=(\tilde Y,Y_5)=(Y_1,\cdots,Y_4,Y_5)$ and $R=|\tilde X|=\sqrt{X_1^2+\cdots+X_4^2}$.

First,  $II_2$  can be estimated as follows.
\begin{eqnarray}\label{Feb22-5}
&&|II_2|\le\int_S^T\int_{R^5}
(T-\tau)^{-\frac52}e^{-C\frac{|X-Y|^2}{T-\tau}}|g^{(k)}(Y,\tau)|
dYd\tau\nonumber\\
&&=\int_S^T [\int_{\{\frac{|X_5-Y_5|^2}{T-\tau}\le C_0,\frac{|\tilde
X-\tilde Y|^2}{T-\tau}\le
C_0\}}+\int_{\{\frac{|X_5-Y_5|^2}{T-\tau}\ge C_0,\frac{|\tilde
X-\tilde Y|^2}{T-\tau}\le
C_0\}}\nonumber\\
&&+\int_{\{\frac{|X_5-Y_5|^2}{T-\tau}\le C_0,\frac{|\tilde X-\tilde
Y|^2}{T-\tau}\ge C_0\}}+\int_{\{\frac{|X_5-Y_5|^2}{T-\tau}\ge
C_0,\frac{|\tilde X-\tilde Y|^2}{T-\tau}\ge
C_0\}}](T-\tau)^{-\frac52}e^{-C\frac{|X-Y|^2}{T-\tau}}|g^{(k)}(Y,\tau)|
dYd\tau\nonumber\\
&&\equiv J_1+J_2+J_3+J_4.
\end{eqnarray}
For $0<\alpha<1$, direct estimates lead to
\begin{eqnarray}\label{Feb22-6}
&&J_1\le \int_S^T \int_{\{|Y|\le
2\sqrt{C_0(T-S)}+2|X|\}}(T-\tau)^{-\frac52}e^{-C\frac{|X-Y|^2}{T-\tau}}\frac{|u_\theta^{(k)}\partial_Zu_\theta^{(k)}|}{R^2}
dYd\tau\nonumber\\
&&\le \int_S^T\int_{\{|Y|\le
2\sqrt{C_0(T-S)}+2|X|\}}(T-\tau)^{-\frac52}e^{-C\frac{|X-Y|^2}{T-\tau}}\frac{|\partial_Zu_\theta^{(k)}|^\alpha}{R^\alpha}
\frac{|\partial_Zu_\theta^{(k)}|^{1-\alpha}}{(R^2|Y|)^{\frac{1-\alpha}{2}}}|Y|^{\frac{1-\alpha}{2}}
\frac{|u_\theta^{(k)}|}{R}dYd\tau\nonumber\\
&&\le(
2\sqrt{C_0(T-S)}+2|X|)^{\frac{1-\alpha}{2}}(\int_S^T\int_{\{|Y|\le
2\sqrt{C_0(T-S)}+2|X|\}} \frac{|\partial_Zu_\theta^{(k)}|^2}{R^2|Y|}
dYd\tau)^{\frac{1-\alpha}{2}}\nonumber\\
&&\times(\int_S^T\int_{\{|Y|\le 2\sqrt{C_0(T-S)}+2|X|\}}
(T-\tau)^{-\frac{5}{1+\alpha}}e^{-C\frac{|X-Y|^2}{T-\tau}}dYd\tau)^{\frac{1+\alpha}{2}}\nonumber\\
&&\le F(C_0,T-S, |X|,\alpha,k)(\int_S^T
(T-\tau)^{-\frac{5}{1+\alpha}+\frac52}d\tau)^{\frac{1+\alpha}{2}}\nonumber\\
&&\le F(C_0,T-S, |X|,\alpha,k),
\end{eqnarray}
where
$$
F(C_0,T-S, |X|,\alpha,k)=F_1(C_0,T-S, |X|,\alpha)
(\int_S^T\int_{\{|Y|\le 2\sqrt{C_0(T-S)}+2|X|\}}
\frac{|\partial_Zu_\theta^{(k)}|^2}{R^2|Y|}
dYd\tau)^{\frac{1-\alpha}{2}} $$  and $F_1(C_0,T-S, |X|,\alpha)$ is a
constant depending on $ C_0,T-S, |X|$ and $\alpha$. Note that the angular component of the velocity in \eqref{Jan2} is
$$
u_\theta^{(k)}(Y,T)=\frac{1}{M_k}u_\theta(\frac{R}{M_k},
z_k+\frac{Y_5}{M_k},t_k+\frac{T}{M_k^2}),
$$
where $z_k=x_{3k}$. Letting
$$
y_1=\frac{Y_1}{M_k}, \cdots,
y_4=\frac{Y_4}{M_k},z=z_k+\frac{Y_5}{M_k}, t=t_k+\frac{\tau}{M_k^2},
r=\sqrt{y_1^2+\cdots+y_4^2},
$$
one can get that
\begin{eqnarray}\label{Nov30-2}
&& \int_S^T\int_{\{|Y|\le 2\sqrt{C_0(T-S)}+2|X|\}}
\frac{|\partial_Zu_\theta^{(k)}|^2}{R^2|Y|} dYd\tau\nonumber\\
&&=\int_{t_k+\frac{S}{M_k^2}}^{t_k+\frac{T}{M_k^2}}\int_{\{\sqrt{r^2+(z-z_k)^2}\le
\frac{2\sqrt{C_0(T-S)}+2|X|}{M_k}\}}
\frac{|\partial_zu_\theta|^2}{\sqrt{r^2+(z-z_k)^2}} rdrdzdt
\end{eqnarray}
Note that $|z_k|=|x_{k3}|\le \bar R$ is bounded. There exists a
subsequence of $\{(z_k, t_k)\}$, still denoted by itself, and  $\bar z\in [-\bar R, \bar R]$ such that
$t_k\to 0, z_k\to \bar{z}$ as $k\to \infty$. By the assumption \eqref{Oct26-2}, for any
$\varepsilon>0$, there exists a $\delta>0$ such that
\begin{eqnarray}\label{Nov-2-0}
 \int_{-2\delta}^{0}\int_0^{2\delta}\int_{\bar z-2\delta}^{\bar z+2\delta}\frac{|\partial_zu_\theta|^2}{\sqrt{r^2+(z-\bar
z)^2}} r drdzdt \le \varepsilon.
\end{eqnarray}
Using \eqref{Oct26-2} again leads to
\begin{eqnarray}\label{Nov-30-0}
 \int_{-2\delta}^{0}\int_0^{2\delta}\int_{\bar z-2\delta}^{\bar z+2\delta}\frac{|\partial_zu_\theta|^2}{\sqrt{r^2+(z-\bar
z_k)^2}} r drdzdt \le C.
\end{eqnarray}
Taking $k>K$ large enough such that
$(t_k+\frac{S}{M_k},t_k+\frac{T}{M_k})\subset (-\delta, 0)$ and
$\{(r,z)\in (0,\infty)\times(-\infty,\infty)|\sqrt{r^2+(z-\bar z)^2}\le
\frac{2\sqrt{C_0(T-S)}+2|X|}{M_k}\}\subset (0,\delta)\times (\bar z-\delta,\bar z+\delta)$.  Let
$$
h_k(r,z,t)=\frac{|\partial_zu_\theta|^2}{\sqrt{r^2+(z-z_k)^2}},
h(r,z,t)=\frac{|\partial_zu_\theta|^2}{\sqrt{r^2+(z-\bar z)^2}},
$$
where $(r,z)\in (0,2\delta)\times(\bar z-2\delta,\bar z+2\delta)$.

For any fixed $0<\delta_0<\delta$, we choose
$0\le\varphi(r,z,t)\le 1$ to be a  smooth function defined in $(0,\infty)\times(-\infty,\infty)\times
(-1,0)$ satisfying $\varphi(r,z,t)\equiv 1$ if $(r,z,t)\in (\delta_0,\delta)\times(\bar z-\delta,\bar z+\delta)\times(-\delta,-\delta_0)$ and
$\varphi(x,t)\equiv 0$ if $(r,z,t)\not\in Q_{\delta,\delta_0}\equiv (\frac{\delta_0}{2},2\delta)\times(\bar z-2\delta,\bar z+2\delta)\times (-2\delta,-\frac{\delta_0}{2})$.

 Then, there exists a subsequence (still denoted by itself) such that
\begin{eqnarray*}
&&h_k\rightharpoonup \tilde h\quad {\rm in}\quad \mathcal{M}(Q_{\delta,\delta_0})\quad\quad (by\ \eqref{Nov-30-0})\\[3mm]
&&h_k\to h\quad {\rm a.e.\ on} \ Q_{\delta,\delta_0},
\end{eqnarray*}
as $k\to \infty$, where $\tilde h\in \mathcal{M}(Q_{\delta,\delta_0})$ which is the finite Radon measure space restricted on $Q_{\delta,\delta_0}$. In particular, it concludes that $\tilde h=h$ and
\begin{eqnarray}\label{Nov-2-1}
\int_{Q_{\delta,\delta_0}}
h_k\varphi
r drdzdt\to \int_{Q_{\delta,\delta_0}}
h\varphi r drdzdt,
\end{eqnarray}
as $k\to \infty$. That is
\begin{eqnarray}\label{Nov-30-1}
\int_{-2\delta}^{-\frac{\delta_0}{2}}\int_{\frac{\delta_0}{2}}^{2\delta}\int_{\bar z-2\delta}^{\bar z+2\delta}
h_k\varphi
r drdzdt\to \int_{-2\delta}^{-\frac{\delta_0}{2}}\int_{\frac{\delta_0}{2}}^{2\delta}\int_{\bar z-2\delta}^{\bar z+2\delta}
h\varphi r drdzdt,
\end{eqnarray}
as $k\to \infty$.
Using \eqref{Nov-2-0} and \eqref{Nov-30-1}, one
obtains, for any $0<\delta_0<\delta$, that
\begin{eqnarray*}
&&\limsup_{k\to\infty}\int_{-\delta}^{-\delta_0}\int_{\delta_0}^{\delta}\int_{\bar z-\delta}^{\bar z+\delta}
h_k
r drdzdt\le \limsup_{k\to\infty}\int_{-2\delta}^{-\frac{\delta_0}{2}}\int_{\frac{\delta_0}{2}}^{2\delta}\int_{\bar z-2\delta}^{\bar z+2\delta}
h_k\varphi
r drdzdt\\
&&=\int_{-2\delta}^{-\frac{\delta_0}{2}}\int_{\frac{\delta_0}{2}}^{2\delta}\int_{\bar z-2\delta}^{\bar z+2\delta}
h\varphi
r drdzdt \le \int_{-2\delta}^{0}\int_0^{2\delta}\int_{\bar z-2\delta}^{\bar z+2\delta} h r drdzdt\le \varepsilon.
\end{eqnarray*}
Due to the arbitrariness of  $\delta_0$ and $\varepsilon>0$, and thanks to \eqref{Feb22-6} and \eqref{Nov-30-1}, one obtains that
\begin{eqnarray}\label{N-6-1}
|J_1|\le F_1(C_0,T-S, |X|,\alpha)\int_{-\delta}^0\int_0^{\delta}\int_{\bar z-\delta}^{\bar z+\delta}
h_k
r drdzdt\to 0, \quad k\to \infty.
\end{eqnarray}

Now we continue estimating $J_2-J_4$.
\begin{eqnarray*}
&&J_2=\int_S^T \int_{\{\frac{|X_5-Y_5|^2}{T-\tau}\ge
C_0,\frac{|\tilde X-\tilde Y|^2}{T-\tau}\le
C_0\}}(T-\tau)^{-\frac52}e^{-C\frac{|X-Y|^2}{T-\tau}}|g^{(k)}(Y,\tau)|
dYd\tau\\
&&\le \int_S^T (T-\tau)^{-\frac52}
d\tau\int_{\{\frac{|X_5-Y_5|^2}{T-\tau}\ge
C_0\}}e^{-C\frac{|X_5-Y_5|^2}{T-\tau}}dY_5\int_{\{\frac{|\tilde X-\tilde
Y|^2}{T-\tau}\le C_0\}}e^{-C\frac{|\tilde X-\tilde Y|^2}{T-\tau}}
\frac{|u_\theta^{(k)}\partial_Zu_\theta^{(k)}|}{R^2} d\tilde Y\\
&&\le \int_S^T (T-\tau)^{-2} d\tau \int_{|\xi|\ge C_0} e^{-C|\xi|^2}
d\xi\int_{\{|\tilde Y|\le \sqrt{C_0(T-S)}+|\tilde
X|\}}e^{-C\frac{|\tilde X-\tilde Y|^2}{T-\tau}}
\frac{|u_\theta^{(k)}\partial_Zu_\theta^{(k)}|}{R^2} d\tilde Y\\
 &&\le C\int_{|\xi|\ge C_0}e^{-C|\xi|^2}
d\xi (T-S),
\end{eqnarray*}
where one has used the transformation $\xi=\frac{X_5-Y_5}{T-\tau}$.
\begin{eqnarray*}
&&J_3=\int_S^T \int_{\{\frac{|X_5-Y_5|^2}{T-\tau}\le
C_0,\frac{|\tilde X-\tilde Y|^2}{T-\tau}\ge
C_0\}}(T-\tau)^{-\frac52}e^{-C\frac{|X-Y|^2}{T-\tau}}|g^{(k)}(Y,\tau)|
dYd\tau\\
&&\le \int_S^T (T-\tau)^{-\frac52}
d\tau\int_{\frac{|X_5-Y_5|^2}{T-\tau}\le
C_0}e^{-C\frac{|X_5-Y_5|^2}{T-\tau}}dY_5\int_{\{\frac{|\tilde
X-\tilde Y|^2}{T-\tau}\ge C_0\}}e^{-C\frac{|\tilde X-\tilde Y|^2}{T-\tau}}|g^{(k)}(Y,\tau)| dY\\
&&\le C(T-S)\int_{|\tilde\xi|\ge C_0} e^{-C|\tilde\xi|^2} d\tilde\xi,
\end{eqnarray*}
where $\tilde\xi=\frac{\tilde X-\tilde Y}{T-\tau}$. Similarly, one can estimate
\begin{eqnarray*}
&&J_4=\int_S^T \int_{\{\frac{|X_5-Y_5|^2}{T-\tau}\ge
C_0,\frac{|\tilde X-\tilde Y|^2}{T-\tau}\ge
C_0\}}(T-\tau)^{-\frac52}e^{-C\frac{|X-Y|^2}{T-\tau}}|g^{(k)}(Y,\tau)|
dYd\tau\\
&&\le C(T-S)\int_{|\xi|\ge C_0} e^{-C|\xi|^2}
d\xi\int_{|\tilde\xi|\ge C_0} e^{-C|\tilde\xi|^2} d\tilde\xi.
\end{eqnarray*}
Putting estimates of $J_1-J_4$ into \eqref{Feb22-6} gives
\begin{eqnarray}\label{Nov30-3}
|II_2|&&\le F(C_0,T-S,
|X|,\alpha,k)\nonumber\\[3mm]
&&+C[(T-S)+1](\int_{|\xi|\ge C_0} e^{-C|\xi|^2} d\xi
+\int_{|\tilde\xi|\ge C_0} e^{-C|\tilde\xi|^2} d\tilde\xi).
\end{eqnarray}
The term $II_1$ can be  treated as for  $J_1$ so that
\begin{eqnarray}\label{Nov30-51}
|II_1|\le F(C_0,T-S,
|X|,\alpha,k)\to 0, k\to\infty,
\end{eqnarray}
for $C_0\ge 2$. Now taking $T-S=C_0\ge 2$ and using \eqref{22},\eqref{24},\eqref{Nov30-3} and \eqref{Nov30-51}, one gets
\begin{eqnarray}\label{N-6-2}
 |f^{(k)}(X,T)|&&\le C[C_0^{-\frac{3}{10}}+F(C_0,C_0,
|X|,\alpha,k)\nonumber\\
&&+C_0(\int_{|\xi|\ge C_0} e^{-C|\xi|^2} d\xi+\int_{|\tilde\xi|\ge C_0} e^{-C|\tilde\xi|^2} d\tilde\xi)].
\end{eqnarray}
Passing to the limit
$k\to \infty$ first and then letting $C_0\to \infty$ in
\eqref{N-6-2}, one obtains that, for any $X\in \mathbb{R}^3, T\in (-\infty,0)$,
\begin{eqnarray*}\label{26+}
|f^{(k)}(X,T)| \to 0, \ k\to \infty,
\end{eqnarray*}
and hence
\begin{eqnarray}\label{N-6-3}
|\omega_\theta^{(k)}(X,T)| \to 0, \ k\to\infty,
\end{eqnarray}
which implies that $\bar\omega_\theta(R,Z,T)=0$. Since $curl_X b^{(k)}(X,T)=\omega_\theta^{(k)}e_\theta$ and
$div_X b^{(k)}(X,T)=0$ so $\bar
b(X,T)$ is a harmonic and bounded function defined on $\mathbb{R}^3\times (-\infty,0)$. Since $\bar u_R(0,Z,T)=0$, thus there exists a continuous and bounded function $s(T): (-\infty,0]\to \mathbb{R}$ such that $\bar b(X,T)=s(T)e_Z$. Using \eqref{Oct26-2}, one obtains that $\bar b(X,T)=0$. To prove that $\bar u_\theta=0$, one needs the following lemma due to \cite{KNSS}.

{\bf Lemma 3.3} Let $u(x,t)$ be a bounded weak solution of the
Navier-Stokes equations in $\mathbb{R}^3\times(-\infty,0)$. Assume that $u(x,t)$
is axisymmetric and satisfies
$$
|u(x,t)\le \frac{C}{\sqrt{x_1^2+x_2^2}} \quad {\rm in}\quad
\mathbb{R}^3\times (-\infty,0).
$$
Then $u=0$ in $\mathbb{R}^3\times (-\infty,0)$.

Applying the facts that $\bar b(X,T)=0$ and $ |R\bar u_\theta|\le C$, the proof of Theorem 3.1 follows from Lemma 3.3.
\hspace{2mm}

\section{Further Remarks}
To rule out the possible singularity of the solution, one should make further investigations. As it is known, Liouville type of theorems play an important role in blow-up approach to study the global regularity of the three-dimensional Navier-Stokes equations. How to apply Theorems 2.1 and Theorem 3.1 to study the global regularity of the three-dimensional axisymmetric Navier-Stokes equations would be interesting.

As mentioned in the proof of Theorems 2.1 and  3.1, if the solution is not bounded, then there exist $t_k\nearrow 0$ as $k\to\infty$ such that
$H(t_k)=h(t_k)$, where $h(t)$ and $H(t)$ are defined as \eqref{P1}.  Denote $N_k=H(t_k)$. Then one can choose a sequence
of numbers $\gamma_k\searrow 1$ as $k\to\infty$ and $x_k\in \mathbb{R}^3$
such that $M_k=|u(x_k,t_k)|\ge N_k/\gamma_k, k=1,2\cdots$,
satisfying $M_k\to \infty$ as $k\to\infty$. Denote $r_k=\sqrt{(x_{k1})^2+(x_{k2})^2}$.

We consider the following two cases.

{\bf Case I.}  $\{r_kM_k\} (k=1,2,\cdots)$ is uniformly bounded.

In this case,  there exists a constant $C>0$ such that
\begin{eqnarray}\label{P3}
r_kM_k\le  C.
\end{eqnarray}
We will use the scaling as in \eqref{P5}.  In view of (\ref{P3}), there exists a point
$(x_{\infty_1}, x_{\infty2})\in \mathbb{R}^2$  such that, up to a
subsequence, $(M_kx_{k1},M_kx_{k2})\to (x_{\infty_1}, x_{\infty2})$
as $k\to\infty$. Here $\sqrt{(x_{\infty_1})^2+(x_{\infty_2})^2}\le
C<\infty$. It follows from (\ref{P8}) and (\ref{27}) that $\bar
u(x_{\infty_1},x_{\infty_2},0)=1$.

Theorem 2.1 implies that $\bar u(X,T)=(0,0, s(T))$, where $s(T)$ is a bounded and continuous function defined in $(-\infty,0]$.  Hence, in this case, to rule out the singularity of the solution, it suffices to prove that $s(T)=0$. Since otherwise, one will obtain a contradiction to the fact that $\bar
u(x_{\infty_1},x_{\infty_2},0)=1$.

Meanwhile, under assumptions of Theorem 3.1, one has $\bar u(X,T)=0$ which is a  contradiction to the fact that $\bar
u(x_{\infty_1},x_{\infty_2},0)=1$ and hence the singularity can be ruled out.

{\bf Case II.}  $r_kM_k (k=1,2,\cdots)$ is not uniformly bounded.

In this case,  one can rescale the solution as
\begin{eqnarray}\label{P1+}
u^{(k)}(X,T)=\frac{1}{M_k}u(x_k+\frac{X}{M_k}, t_k+\frac{T}{M^2_k}).
\end{eqnarray}
 Then, similar to Case I, $u^{(k)}(X,T) (k=1,2,\cdots)$ are  smooth solutions of the 3D Navier-Stokes equations, which are defined in
$\mathbb{R}^3\times (A_k,B_k)$ with
\begin{eqnarray}\label{P6-}
A_k=-M_k^2-M_k^2t_k, \ B_k=-M_k^2t_k.
\end{eqnarray}
Note that $B_k=-M_k^2t_k\ge (\frac{N_k}{\gamma_k})^2(-t_k)\ge
\frac{\varepsilon_1^2}{\gamma_k^2}$ for some $\varepsilon_1>0$. Moreover, it is clear that
\begin{eqnarray}\label{P7-}
|u^{(k)}(X,T)|\le \gamma_k, X\in \mathbb{R}^3, T\in (A_k,0],
\end{eqnarray}
and
\begin{eqnarray}\label{P8-}
 |u^{(k)}(0,0)|=1.
\end{eqnarray}

In this case, there exists a subsequence of $\{M_k\}$ (still denoted
by itself) such that $r_kM_k\to \infty$ as $k\to\infty$. Due to the
axis symmetry of $u$, $x_k$ can be chosen so that $\theta(x_k)\to
\theta_\infty$ for some $\theta_\infty\in [0,2\pi]$. Then there exists an unit vector $\nu=(\nu_1,\nu_2,0)$ such that
$e_r(x_k)\to \nu$ and $e_\theta(x_k)\to \nu^\perp=(-\nu_2,\nu_1,0)$. Moreover, it holds that
$$
 x_k+\frac{X}{ M_k}\in B(x_k, \frac{r_k}{\sqrt{r_kM_k}})\ \ {\rm for}\ \ X\in B(0,\sqrt{r_k M_k}),
$$
and
$$
t_k-(\frac{r_k}{\sqrt{r_k M_k}})^2<t_k+\frac{T}{ M_k^2}\le t_k<0 \ \
{\rm for }\ \ -M_kr_k<T\le 0.
$$
By (\ref{C10-1}),
$$
|u_\theta(y,t)|\le \frac{C}{r_k} \ {\rm for}\ \ y\in B( x_k,
\frac{r_k}{2}), t<0,
$$
which implies that
\begin{eqnarray}\label{33-}
|u^{(k)}(X,T)e_\theta(x_k+\frac{X}{M_k})|=\frac{1}{M_k}|u_\theta(
x_k+\frac{X}{M_k}, t_k+\frac{T}{M_k^2})|\le \frac{C}{M_kr_k}
\end{eqnarray}
for $(X,T)\in B(0,\sqrt{r_kM_k})\times (-r_kM_k,0]$.

Since the flow is axisymmetric, thus,  on
$B(0,\sqrt{r_kM_k})\times (-r_kM_k,0]$, $e_R(x_k+\frac{X}{M_k})\to
\nu$ and $e_\theta(x_k+\frac{X}{ M_k})\to \nu^\perp$ as
$k\to\infty$. Moreover, for each $k$, $u^{(k)}$ is a bounded and smooth solution to the 3D Navier-Stokes equations. There exists a subsequence of $u^{(k)}$ (still denoted by itself) and a bounded
ancient solution $\tilde u(X,T)$ to the 3D Navier-Stokes equations on
$\mathbb{R}^3\times (-\infty,0]$, such that
\begin{eqnarray}
&&u^{(k)}(X,T)=\frac{1}{M_k}u_R(x_k+\frac{X}{
M_k}, t_k+\frac{T}{M_k^2})e_R(x_k+\frac{X}{M_k})\nonumber\\[3mm]
&&+\frac{1}{M_k}u_\theta(x_k+\frac{X}{M_k}, t_k+\frac{T}{M_k^2})e_\theta(x_k+\frac{X}{M_k})\nonumber\\
&&+\frac{1}{M_k}u_Z(x_k+\frac{X}{M_k},\frac{T}{M_k^2})e_Z(x_k+\frac{X}{M_k})\nonumber\\[3mm]
&&\rightarrow \tilde u(X,T)=\tilde u_R\nu+\tilde u_\theta\nu^\perp+\tilde u_Ze_Z\nonumber
\end{eqnarray}
in $C(\bar Q)$ for any compact subset $Q$ of $\mathbb{R}^3\times
(-\infty,0]$.

In view of \eqref{P8-}, one has
\begin{eqnarray}\label{Oct29}
|\tilde u(0,0)|=1.
\end{eqnarray}

 In this case,  Lei and Zhang  \cite{LZ} obtained a Liouville type of theorem which can be stated as follows.

 \hspace{3mm}

 {\bf Proposition 4.1} If $r_kM_k (k=1,2,\cdots)$ is not uniformly bounded, then $\tilde u(X,T)=(0, s_1(T), s_2(T))$, where $s_1(T)$ and $s_2(T)$ are bounded and continuous functions depending on $T\in (-\infty,0]$.

\hspace{3mm}
The proof is referred to \cite{LZ} and we give a sketch of the proof here.

\hspace{3mm}

 {\bf Sketch of Proof}. \eqref{33-} implies  $\tilde u(X,T)\cdot \nu^\perp=0$. Hence
\begin{eqnarray}\label{34-}
\tilde u(X,T)=\tilde u_R(X,T)\nu+\tilde u_Z(X,T)e_Z.
\end{eqnarray}

On the other hand, for $(y,s)\in B(x_k, \frac{r_k}{\sqrt{r_k
M_k}})\times [t_k-(\frac{r_k}{\sqrt{r_kM_k}})^2, t_k]$, one has that
\begin{eqnarray*}
&&\frac{1}{M_k}[u_r(y,s)e_\theta(y)-u_\theta(y,s)e_r(y)]\\[3mm]
&&=\frac{1}{M_k}\partial_\theta[u_r(y,s)e_r(y)+u_\theta(y,s)e_\theta(y)]\\[3mm]
&&=\frac{1}{M_k}\partial_\theta[u_r(y,s)e_r(y)+u_\theta(y,s)e_\theta(y)+u_z(y,s)e_z(y)]\\[3mm]
&&=\partial_\theta[u^{(k)}(M_k(y-x_k), M_k^2(s-t_k))]\\[3mm]
&&=M_k(\partial_\theta y\cdot\nabla)u^{(k)}(M_k(y-x_k),M_k^2(s-t_k))\\[3mm]
&&=M_k|y|(e_\theta(y)\cdot\nabla)u^{(k)}(M_k(y-x_k), M_k^2(s-t_k)),
\end{eqnarray*}
which shows that
\begin{eqnarray}
&&\frac{1}{M_k}[u_r(x_k+\frac{X}{M_k}, t_k+\frac{T}{M_k^2})e_\theta(x_k+\frac{X}{M_k})\nonumber\\[3mm]
&&-u_\theta(x_k+\frac{X}{M_k}, t_k+\frac{T}{M_k^2})e_R(x_k+\frac{X}{M_k}))]\nonumber\\[3mm]
&&=M_k|x_k+\frac{X}{M_k}|(e_\theta(x_k+\frac{X}{M_k})\cdot\nabla)u^{(k)}(X,T)\label{35-}
\end{eqnarray}
for $(X,T)\in B(0,\sqrt{r_kM_k})\times (-r_k M_k,0]$. Since $r_k
M_k\to \infty$, so $M_k|x_k+\frac{X}{M_k}|\to \infty$ for  any fixed
$X\in B(0,\sqrt{r_k M_k})$. But the left hand side of (\ref{35-}) is
bounded. Hence, letting $k\to \infty$, one gets that
\begin{eqnarray}\label{36}
(\nu^\perp\cdot\nabla)\tilde u(X,T)=0.
\end{eqnarray}

Note that the Navier-Stokes equations are invariant under rotation.
Without loss of generality, we set $\nu=e_1, \nu^\perp=e_2$.
Consequently, the limit function
$$
\tilde u(X,T)=\tilde u_R(X_1,Z,T)e_1+\tilde u_Z(X_1,Z,T)e_Z,
$$
is a bounded ancient  solution to the 2D Navier-Stokes
equations. It follows from Theorem 5.1 and Remark 6.1 in \cite{KNSS} that $\nabla \tilde u_R=\nabla \tilde u_Z=0$. Hence, $\tilde u_R$ and $\tilde u_Z$ are bounded and continuous functions depending only on time variable $T\in(-\infty,0]$, denoted by $s_1(T)$ and $s_2(T)$ respectively. The proof of the proposition is finished.

\hspace{3mm}

We finally remark that in  case II, to role out the singularity of the solution, it suffices to prove that $s_1(T)=s_2(T)=0$. Since otherwise, one can obtain a contradiction to \eqref{Oct29}.

\end{document}